\documentclass[11pt]{article}

\usepackage{amsmath,amsthm,epsf,graphicx,listings}

\newcommand{\EQ}[2]{\begin{equation}{{#2}\label{#1}} \end{equation}}
                        \newcounter{theorem}
\makeatletter \@addtoreset{theorem}{section} \makeatother

\renewcommand{\thetheorem}{\arabic{section}.\arabic{theorem}}

\newenvironment{thm}[2]{\begin{sloppypar}\refstepcounter{theorem}%
                        {\bf #1 \thetheorem.}\label{#2}\em{}}%
                        {\end{sloppypar}}
                        \newcommand{\theo}[3]{\begin{thm}{#1}{#2}
#3\end{thm}}

                        \newcommand{\N}{{\rm I}\!{\rm N}}
             
                         \newcommand{\Z}{\makebox[0.04cm][l]{\sf Z}\!{\sf Z}}

                        \newcommand{\R}{{\rm I}\!{\rm R}}

\def\endofproof{\vbox{\hrule\hbox{\vrule\hbox to 6truept
                                      {\vbox to 6truept
{\vfil}\hfil}\vrule}\hrule}}%
 \numberwithin{equation}{section}

                        \renewcommand{\author}{Aurelian Bejancu and Simon Hubbert\footnote{Corresponding author. E-mail: s.hubbert@bbk.ac.uk} }

                        \renewcommand{\title}{A  study of
                        the uniform accuracy of univariate thin plate spline interpolation}  

\begin{document}
                        \setcounter{page}{1}
                        \vskip1cm
                        \begin{center}
                        \LARGE{\bf \title}
                        \\[0.7cm]

            \large{\author}
                        \end{center}
                        \vspace{0.3cm}

                        \begin{abstract}

The usual power function error estimates do not capture the true
order of uniform accuracy for thin plate spline interpolation to
smooth data functions in one variable. In this paper we propose a
new type of power function and we  show, through numerical
experiments, that the error estimate based upon it does match the
expected order.
 We also study the relationship between the new power function
 and the Peano kernel for univariate thin plate spline interpolation.


                        \end{abstract}

\setcounter{section}{0} 

\section{Introduction}

For each $\gamma>0$, define the basis function $\phi_{\gamma}:\R
\rightarrow \R$ by
\[
\phi_{\gamma}\left(  x\right)  =\left\{
\begin{array}
[c]{lc}%
|x|^{\gamma}, & \text{if }\gamma\not \in 2\N,\\
|x|^{\gamma}\log |x|, & \text{if }\gamma\in2\N.
\end{array}
\right.
\]
Let $m_{\gamma}= \lfloor \gamma/2 \rfloor$ be the integer part of
$\gamma/2$. For any integer $n \geq m_{\gamma}$ and any set of
values $\left\{ f_{0},\ldots,f_{n}\right\}  $ of a target function
$f$ prescribed at the set of equi-spaced knots $\left\{
0,h,2h,\ldots,1\right\}$, where $h=1/n$, Micchelli's theory
\cite{Mic} of conditionally positive definite radial basis functions
guarantees the existence of a unique function $s_{h,\gamma}$ of the
form
\EQ{rbfint}{ s_{h,\gamma}\left(  x\right)
=\sum_{k=0}^{n}a_{k}\phi_{\gamma} \left(  x-hk \right) +
\sum_{l=0}^{m_{\gamma}}b_{l}x^{l},\quad x\in\R,}
that satisfies the interpolation conditions
\EQ{interpcond}{ s_{h,\gamma}\left(  hi\right)  =f_{i},\quad
i=0,1,\ldots,n,}
as well as the `side' conditions
\EQ{sideconds}{ \sum_{k=0}^{n}a_{k}k^{l}=0,\quad
l=0,1,\ldots,m_{\gamma}. }

In this paper we focus on the special case $\gamma=2$ corresponding
to the thin plate spline (TPS) basis function $\phi_{2}\left( x
\right) = |x|^{2}\log |x|$ and we investigate the rate at which the
interpolant (\ref{rbfint}) converges to $f$ uniformly over $[0,1]$,
as $h\to 0$. If $f$ has a Lipschitz continuous third derivative on
$[0,1]$, Bejancu \cite{Bej2} proved that the uniform error over a
fixed compact subset of $(0,1)$ inherits the maximal convergence
rate $O(h^{3})$ obtained by Powell \cite{P92} and Buhmann
\cite{BuCA} for `cardinal' TPS interpolation on the infinite grid $h
\Z$. Due to boundary effects, however, the uniform norm of the error
over the full interval $[0,1]$ decays at the much slower rate
$O(h^{3/2})$, as illustrated numerically in \cite{B99,Pow,hm}.

The usual method for error estimation in radial basis function
interpolation, reviewed in the next section, delivers bounds of the
form
 \[
|f(x)-s_{h,\gamma}\left( x\right)|\le
 c_{f,\gamma}\mathcal{P}_{h,\gamma}(x),\quad  x \in [0,1],
 \]
where $\mathcal{P}_{h,\gamma}$ is the so-called `power function'
associated with $\phi_{\gamma}$, while $f$ belongs to the `native
space' generated by $\phi_{\gamma}$ \cite{WendBook,WS}. It is well
known that theoretical convergence rates based upon bounding
$\mathcal{P}_{h,\gamma}(x)$ uniformly for $x\in [0,1]$ do not match
the actual rates of decay of the error achieved in numerical
experiments if $f$ has sufficiently many continuous derivatives.
This discrepancy was first observed by Powell \cite{P94} for the
bivariate TPS interpolant.

For $\gamma=2$, in section 3, we obtain a new error bound which
employs a `mixed power function' $\mathcal{M}_{h,\mu}$ defined by
means of the basis functions $\phi_2$ and $\phi_{\mu}$, for $\mu \in
(0,4)$. We then perform a numerical study of $\max_{x\in [0,1]}
\mathcal{M}_{h,\mu} (x)$ as $h\rightarrow 0$, which shows that, for
$\mu \in [3,4)$, the mixed power function decays like a constant
multiple of $h^{3/2}$. This matches exactly the previously known
numerical order of uniform convergence of the error $f-s_{h,2}$ on
$[0,1]$, for sufficiently smooth target functions $f$. In section 4
we prove that, for $\mu = 3$ and $x\in [0,1]$, the mixed power
function value $\mathcal{M}_{h,3} (x)$ is, up to a constant factor,
the $L^2$-norm of the Peano kernel of the error functional at $x$.
Moreover, we provide numerical evidence that the smaller $L^1$-norm
of this Peano kernel does not in fact decay faster than the mixed
power function when measured uniformly over $[0,1]$. It is hoped
that these results and the conjectures formulated in the paper will
motivate future work to establish theoretically the uniform
convergence order $O(h^{3/2})$ for univariate TPS interpolation to
sufficiently smooth target functions.

\section{Error estimates via the standard power function}

In this section we review the power function technique to obtain
error estimates for univariate interpolation with the radial basis
function $\phi_{\gamma}$.  A key role in this technique is played by
the generalized or distributional Fourier transform of
$\phi_{\gamma}$.

\theo{Lemma}{ftrbf}{{\rm{\textbf{\cite[section\ 8.3]{WendBook}}}}
For each $\gamma > 0$, the generalized Fourier transform of
$\phi_{\gamma}$ satisfies
\[
\widehat{\phi_{\gamma}}(t) = \frac{A_{\gamma}}{|t|^{1+\gamma}},\quad
t \in \R\setminus \{0\},\]
for some constant $A_{\gamma}$ such that
$(-1)^{m_{\gamma}+1}A_{\gamma}>0$.}

\subsection{The standard power function}

As above, let $m_{\gamma} = \lfloor\gamma/2\rfloor,$ $ n \ge
m_{\gamma}$, and $h = 1/n.$ For each $x \in \R$ which is not in the
knot-set $\{0,h,\ldots,1\}$, Micchelli's theory implies that the
quadratic form
\EQ{dirpow}{
Q_{\gamma,n}(\mathbf{v}):=(-1)^{m_{\gamma}+1}\left(\sum_{j=0}^{n}\sum_{k=0}^{n}v_{j}v_{k}
\phi_{\gamma}(hj-hk)-2\sum_{j=0}^{n}v_{j}\phi_{\gamma}(x-hj)\right),}
is strictly positive whenever the non-zero vector $\mathbf{v} =
(v_{0},\ldots,v_{n})^{T}\in \R^{n+1}$ satisfies
\EQ{cons}{x^{l}=\sum_{j=0}^{n}v_{j}\left(  hj\right) ^{l}  ,\quad \
l=0,1,\ldots,m_{\gamma}.}
Further, for each $j \in \{0,1,\ldots,n\},$ let
$\ell_{j,h}^{(\gamma)}$ be the unique function of the type
(\ref{rbfint})--(\ref{sideconds}) which satisfies the Lagrange
interpolation conditions
\[
\ell_{j,h}^{(\gamma)}\left(  ih\right)  =\delta_{ij},\quad
i=0,1,\ldots,n.
\]
Then we have the Lagrange representation formula for the interpolant
(\ref{rbfint}):
\begin{equation}
s_{h,\gamma}\left(  x\right)
=\sum_{j=0}^{n}f_{j}\ell_{j,h}^{(\gamma)}\left( x\right)  ,\quad
x\in\R,\label{eq:Lagrange}%
\end{equation}
as well as the reproduction formula
\EQ{repro}{ x^{l}=\sum_{j=0}^{n}\left(  hj\right)
^{l}\ell_{j,h}^{(\gamma)}\left( x\right)  ,\quad x\in\R,\,\,\,
l=0,1,\ldots,m_{\gamma}. }
\theo{Proposition }{variational}{{\rm{\textbf{\cite{WS}}}}\ With the
above notations, for each $x \in \R$, the vector
\[
\mathbf{v}_{x}= \Bigl(\ell_{0,h}^{(\gamma)}\left( x\right),\ldots,
\ell_{n,h}^{(\gamma)}\left( x\right)\Bigr)^{T}\in \R^{n+1}
\]
has the property that it minimizes the quadratic form
{\rm{(\ref{dirpow})}} among all non-zero vectors $\mathbf{v} \in
\R^{n+1}$ that satisfy {\rm{(\ref{cons})}}.}

The minimum value of the quadratic form $Q_{\gamma,n}$ defines the
square of the so-called `power function'
$\mathcal{P}_{h,\gamma}:\R\to [0,\infty)$, namely
\EQ{dirpowdef}{
\mathcal{P}_{h,\gamma}^{2}(x):=Q_{\gamma,n}(\mathbf{v}_{x}).}
Note that $\mathcal{P}_{h,\gamma} (hj) = 0$, $\forall
j\in\{0,1,\ldots,n\}$.

\theo{Proposition}{intrepth}{{\rm{\textbf{\cite{WS}}}} \ For each $x
\in \R$, let
\EQ{theta}{ \Theta_{x,\gamma}\left( t\right)
:=e^{ixt}-\sum\limits_{j=0}^{n}\ell_{j,h}^{(\gamma)}\left( x\right)
e^{ihjt},\quad t\in \R.}
Then we have the absolutely convergent integral representation
\EQ{intrep}{
 \mathcal{P}_{h,\gamma}^{2}(x)
 = \frac{|A_{\gamma}|}{2\pi} \int_{\R}
\frac{ \left\vert\Theta_{x,\gamma}\left(t\right)\right\vert^{2} } {
|t|^{1+\gamma} } dt .} }

\subsection{Error estimates}

For each $\gamma>0$, let $\kappa_{\gamma} = \lceil \gamma/2+1\rceil$
be the least integer that is greater than or equal to $\gamma/2+1$.
In order to obtain error bounds for any target function $f \in
C^{\kappa_{\gamma}}[0,1]$, we construct an extension $f^{*}:\R\to\R$
of $f$ as follows (cf.\ \cite{Bej2}). By the Whitney extension
theorem \cite{Whit}, there exists $\widetilde{f} \in
C^{\kappa_{\gamma}}(\R)$ such that $\widetilde{f}(x)  =f(x)$, for $x
\in [0,1]$. Let $\nu$ be an infinitely differentiable cut-off
function which satisfies $\nu(x)=1$ for $x \in [0,1]$ and $\nu(x)=0$
for sufficiently large $|x|$, and set
\[
f^{*}(x) := \nu(x)\widetilde{f}(x),\quad \forall x \in \R.
\]
Clearly $f^{*} \in C^{\kappa_{\gamma}}(\R)$ is compactly supported
and coincides with $f$ on $[0,1]$. Furthermore, its Fourier
transform $\widehat{f^{\ast}}$, defined as the continuous function
\[
\widehat{f^{\ast}}\left(  t\right)
:=\int_{\R}e^{-ixt}f^{\ast}\left(  x\right)  dx,
\]
satisfies
\EQ{new}{ \left\vert t^{\kappa_{\gamma}}\widehat{f^{\ast}}\left(
t\right) \right\vert
=\left\vert \widehat{\left(  f^{\ast}\right)  ^{\left(  \kappa_{\gamma}\right)  }%
}\left(  t\right)  \right\vert \leq\left\Vert \left( f^{\ast}\right)
^{\left(  \kappa_{\gamma}\right)  }\right\Vert _{L^{1}\left(
\R\right)  },}
for any $t\not =0$. In particular, $\widehat{f^{\ast}}$ is
integrable over $\R$, so $f^{\ast}$ can be recovered via the Fourier
inversion formula
\EQ{finv}{ f^{*}\left(  x\right)
=\frac{1}{2\pi}\int_{\R}e^{ixt}\widehat{f^{*}}\left( t\right)
dt,\quad x\in\R. }

Next, let $f_{j}:=f\left( hj \right)$ in (\ref{interpcond}) for
$j=0,\ldots,n$. Then (\ref{finv}) and (\ref{eq:Lagrange}) imply the
error representation
\EQ{error}{ f\left(  x\right)  -s_{h,\gamma}\left(  x\right)
=\frac{1}{2\pi}\int_{\R } \widehat{f^{*}}\left(
t\right)\Theta_{x,\gamma}\left( t\right) dt,\quad x\in\left[
0,1\right],}
where $\Theta_{x,\gamma}$ is given by (\ref{theta}). Moreover, as a
consequence of (\ref{new}) and the definition of $\kappa_{\gamma}$,
we have
\[
c_{f,\gamma} :=  \Bigl\{\frac{1}{2\pi |A_{\gamma}|}
 \int_{\R}|\widehat{f^{*}}(t)|^{2}|t|^{1+\gamma}dt\Bigr\}^{1/2}<\infty,
\]
i.e., $f^{*}$ belongs to the so-called `native space' generated by
$\phi_{\gamma}$. Using (\ref{intrep}) and the Cauchy-Schwarz
inequality in (\ref{error}), we obtain the error bound
\EQ{stdpow-bnd}{ \left\vert f\left(  x\right)  -s_{h,\gamma}\left(
x\right) \right\vert \leq c_{f,\gamma}\mathcal{P}_{h,\gamma}\left(
x\right) ,\quad x\in [0,1]. }
Further, Wu and Schaback \cite{WS} showed that the variational
characterization of the power function given in
Proposition~\ref{variational} implies
\EQ{aost}{ \max_{x\in\left[  0,1\right] }\mathcal{P}_{h,\gamma}(x)
\leq B_{\gamma}h^{\gamma/2},\quad\text{as }h\rightarrow 0,}
for a constant $B_{\gamma}$ independent of $h$. On the other hand,
Schaback and Wendland \cite{SW} proved that the exponent $\gamma/2$
cannot be increased in the above bound. Thus, the power function
technique leads to the maximal estimate $O(h^{\gamma/2})$ for the
uniform norm of the error over $[0,1]$.

\section{A mixed power function for univariate TPS}

In this section we focus on the TPS basis function $\phi_{2}(x) =
|x|^{2}\log |x|$, i.e.\ $\gamma=2$. According to (\ref{aost}), in
this case the power function satisfies
\[
\max_{x\in\left[ 0,1\right]  }\mathcal{P}_{h,2}\left(  x\right)
=O\left(  h\right),\quad \text{as}\ h\rightarrow 0.
\]
However, numerical experiments \cite{B99,hm,Pow} suggest that the
uniform error \EQ{uner}{\max_{x\in\left[ 0,1\right]  }\left\vert
f\left( x\right) -s_{h,2}\left(  x\right) \right\vert} is of the
magnitude of $h^{3/2}$ for a sufficiently smooth target function
$f$.

To address this discrepancy, we start from the integral
representation (\ref{intrep}):
\EQ{intreptps}{
\mathcal{P}_{h,2}^{2}(x)=\frac{|A_{2}|}{2\pi}\int_{\R
}\frac{\left\vert\Theta_{x,2}\left(
t\right)\right\vert^{2}}{|t|^{3}}dt.}

Note that expression (\ref{theta}) satisfies
 \EQ{theta-behav}{ \left\vert
\Theta_{x,2}\left( t\right) \right\vert =\left\{
\begin{array}
[c]{ll}%
O\left(  \left\vert t\right\vert ^{2}\right)  , & \text{as }t\rightarrow0,\\
O\left(  1\right)  , & \text{as }\left\vert t\right\vert
\rightarrow\infty,
\end{array}
\right. } for each fixed $x$ and $h.$ Indeed, since $m_{2}=1$, by
(\ref{repro}) we have
\EQ{replin}{ 1=\sum_{j=0}^{n}\ell_{j,h}^{(2)}\left(x\right)\quad
{\rm{and}}\quad x=\sum_{j=0}^{n}hj\cdot\ell_{j,h}^{(2)}\left(
x\right),\quad {\rm{for}}\ x\in \R. }
This fact, together with the series expansion of the exponential,
provides the bound (\ref{theta-behav}) for $ t \rightarrow 0$.  The
bound for $\left\vert t\right\vert\rightarrow\infty$ follows from
the triangle inequality.

As a consequence of (\ref{theta-behav}), the integral
(\ref{intreptps}) is still well defined if  $|t|^{3}$ is replaced by
$|t|^{1+\mu}$, for any $\mu \in (0,4)$, $\mu\not= 2$. We may thus
define the \emph{mixed power function}
$\mathcal{M}_{h,\mu}:\R\rightarrow\lbrack0,\infty)$ whose square is
given by
\EQ{mixed-pf}{
\mathcal{M}_{h,\mu}^{2}\left(  x\right) :=\frac{|A_{\mu}|%
}{2\pi}\int_{\R}\frac{\left\vert \Theta_{x,2}\left( t\right)
\right\vert ^{2}}{\left\vert t\right\vert ^{1+\mu}}dt,\quad x\in
\R,\,\,\, \mu\in\left(  0,4\right).}
Under this integral, the Lagrange functions entering in expression
(\ref{theta}) of $\Theta_{x,2}$ (and generated by the TPS basis
function $\phi_{2}$) are combined with the generalized Fourier
transform of the basis function $\phi_{\mu}$ (cf.\ Lemma
\ref{ftrbf}).

We now let $\kappa=\left\lceil \mu/2+1\right\rceil $, $f\in
C^{\kappa}\left[ 0,1\right]$ and, as in subsection~2.2, consider the
compactly supported extension $f^{\ast}\in C^{\kappa}\left(
\R\right)$ of $f$ to the whole real axis. Then the error analysis of
subsection~2.2 recast in terms of the mixed power function implies
\EQ{mixed-bound}{ \left\vert f\left(  x\right)  -s_{h,2}\left(
x\right) \right\vert \leq c_{f,\mu}\mathcal{M}_{h,\mu}\left(
x\right) ,\quad x\in\left[  0,1\right] ,\ \mu\in\left(  0,4\right)
,}
where
 \EQ{mixed-const}{ c_{f,\mu}=\Bigl\{\frac{1}{2\pi|A_{\mu}|}
\int_{\R}\left|\widehat{f^{\ast}}\left(t\right) \right|
^{2}|t|^{1+\mu}dt\Bigr\}^{1/2}<\infty.}
This shows that, for a fixed $\mu\in (0,4)$, estimates of the decay
of the mixed power function $\mathcal{M}_{h,\mu}$ as $h\rightarrow
0$ will deliver error estimates for TPS interpolation. Therefore we
state the following problem:

\theo{Problem}{problem}{Given $\mu \in (0,4)$, $\mu\not= 2$, does
there exist an algebraic decay rate of the mixed power function
uniformly over $[0,1]$, i.e., a largest value $\alpha_{\mu}>0$ such
that
\EQ{ao}{ \max_{x\in\left[  0,1\right] }\mathcal{M}_{h,\mu}\left(
x\right) = O\left( h^{\alpha_{\mu}} \right),\quad {\rm{as }}\ h \to
0?} }

Before embarking on a numerical answer to this problem, a few
remarks are in order. Firstly, note that, due to
(\ref{mixed-const}), the above target function $f\in
C^{\kappa}\left[ 0,1\right]  $ has its compactly supported extension
$f^{\ast}$ in the native space generated by the basis function
$\phi_{\mu}$, rather than the native space generated by the TPS
basis function $\phi_{2}$ as in the standard estimate
(\ref{stdpow-bnd}) for $\gamma=2$. Thus, for a given $\mu \in (0,4)$
$(\mu\not= 2),$ the resulting mixed power function error bound
(\ref{mixed-bound}) is precisely what we would expect if we measured
the TPS interpolation error for target functions in the native space
of $\phi_{\mu};$ this approach is investigated in \cite{Schabn} in
the context of approximation rather than interpolation. In
particular, the mixed power function bound (\ref{mixed-bound})
applies to the smooth ($C^{\infty}$) target functions employed in
the numerical experiments of \cite{B99,hm,Pow}.

Secondly, recall the two equivalent expressions for the standard
power function: the direct form (\ref{dirpowdef}) and the integral
representation (\ref{intrep}). Letting $m=\left\lfloor
\mu/2\right\rfloor $, an application of Theorem~3 from \cite{WS}
shows that the mixed power function can also be expressed as
\EQ{mixedpfdir}{
\begin{aligned}
\mathcal{M}_{h,\mu}^{2}(x)=&(-1)^{m+1}\Biggl(\sum_{j=0}^{n}
\sum_{k=0}^{n}\ell_{j,h}^{(2)}\left( x\right)\ell_{k,h}^{(2)}\left(
x\right)\phi_{\mu}(hj-hk)\\
\, & \\&-2\sum_{j=0}^{n}\ell_{j,h}^{(2)}\left(
x\right)\phi_{\mu}(x-hj)\Biggr),\quad x\in\R.
\end{aligned}}

Thirdly, note that (\ref{replin}) implies that the  TPS Lagrange
functions  $\ell_{j,h}^{(2)}$ satisfy constraint  (\ref{cons}) of
the variational problem from Proposition~\ref{variational} with
$\mu$ in place of $\gamma$. However, the solution to that problem is
provided by the values of the Lagrange functions generated by
interpolation with the basis function $\phi_{\mu}$. As a result, the
bounding technique \cite{WS} that leads to the estimate (\ref{aost})
for the standard power function cannot be applied to obtain
estimates on $\mathcal{M}_{h,\mu}$.

We now turn to a numerical investigation of the behaviour of the
mixed power function. For a fixed parameter $\mu\in (0,4)$,
$\mu\not= 2$, we compute an approximation
$\mathcal{M}_{h,\mu}^{(max)}$ of the left-hand side of (\ref{ao})
for $h =1/n$, starting from $n=128$ and proceeding as follows:
\begin{enumerate}

    \item For the current mesh-size $h$ and each $j\in \{0,1,\ldots,n\}$, express
          the TPS Lagrange function $\ell_{j,h}^{(2)}$ in the form (\ref{rbfint}) and compute
      its coefficients by solving the system (\ref{interpcond})--(\ref{sideconds}),
      where $f_i=\delta_{ij}$, $i\in \{0,1,\ldots,n\}$.

    \item Use (\ref{mixedpfdir}) to evaluate the mixed power function at the set of
      mid-points $\mathcal{X}_{eval,h}=\left\{h/2,3h/2,\ldots,1-h/2\right\}$
          and determine its maximum value
    \[
    \mathcal{M}_{h,\mu}^{(max)}=\max\left\{\mathcal{M}_{h,\mu}(x):x \in
    \mathcal{X}_{eval,h}\right\}.
    \]
    \item Double $n$ and repeat steps 1--2 as long as $n\leq 2048$.

\end{enumerate}
\noindent The results displayed in Tables \ref{o1}--\ref{o4} show
that,
for each chosen $\mu$, the values of $\mathcal{M}_{h,\mu}^{(max)}$ satisfy%
\[
\mathcal{M}_{h,\mu}^{(max)}=c_{\mu}h^{\alpha_{h,\mu}},
\]
where $c_{\mu}$ and $\alpha_{h,\mu}$ are also included in the
tables. On the basis of these numerical results, we are led to the
following conjecture.

\begin{table}
\begin{tabular}{| c|c |c | c|c |c|c |}\hline
& \multicolumn{2}{c|}{$\mu = 1/3$} & \multicolumn{2}{c|}{$\mu =
2/3$}& \multicolumn{2}{c|}{$\mu = 1$}\\ \hline
   $h^{-1}$   & $\mathcal{M}_{h,\mu}^{(max)}$& $\alpha_{h,\mu}$ &
    $\mathcal{M}_{h,\mu}^{(max)}$ & $\alpha_{h,\mu}$  & $\mathcal{M}_{h,\mu}^{(max)}$ & $\alpha_{h,\mu}$ \\  \hline
  128 &  4.774E-01& 0.167      &   1.768E-01& 0.333   & 6.342E-02 & 0.500\\
  256 &  4.253E-01& 0.167      &   1.404E-01& 0.333   & 4.485E-02 & 0.500\\
  512 &  3.789E-01 & 0.167  & 1.114E-01& 0.333  & 3.171E-02 & 0.500\\
  1024 & 3.376E-01& 0.167   & 8.842E-02& 0.333  &2.242E-02 &0.500\\
  2048 & 3.007E-01 &0.167   & 7.018E-02& 0.333  &1.586E-02& 0.500\\
  \hline
  $c_{\mu}$ & \multicolumn{2}{c|}{1.072} & \multicolumn{2}{c|}{0.8912} & \multicolumn{2}{c|}{0.7175}\\
  \hline
\end{tabular}
  \centering
  \caption{Decay of the mixed power function for $\mu
   \in (0,1]$}\label{o1}
\end{table}

\begin{table}
\begin{tabular}{| c|c |c | c|c|}\hline
& \multicolumn{2}{c|}{$\mu = 4/3$} & \multicolumn{2}{c|}{$\mu =
5/3$}\\ \hline
   $h^{-1}$   & $\mathcal{M}_{h,\mu}^{(max)}$& $\alpha_{h,\mu}$ &
    $\mathcal{M}_{h,\mu}^{(max)}$& $\alpha_{h,\mu}$  \\  \hline
    128 &  2.143E-02 & 0.667      &   6.258E-03& 0.833   \\
  256 &  1.350E-02& 0.667      &   3.512E-03& 0.833   \\
  512 &  8.503E-03 & 0.667  & 1.971E-03& 0.833  \\
  1024 & 5.356E-03& 0.667   & 1.106E-03& 0.833  \\
  2048 & 3.374E-03 &0.667   & 6.208E-04& 0.833  \\
  \hline
  $c_{\mu}$ & \multicolumn{2}{c|}{0.5442} & \multicolumn{2}{c|}{0.3568} \\
  \hline
\end{tabular}
  \centering
  \caption{Decay of the mixed power function for $\mu
   \in (1,2)$}\label{o2}
\end{table}

\begin{table}
\begin{tabular}{| c|c |c | c|c |c|c |}\hline
& \multicolumn{2}{c|}{$\mu = 7/3$} & \multicolumn{2}{c|}{$\mu =
8/3$}& \multicolumn{2}{c|}{$\mu = 3$}\\ \hline
   $h^{-1}$   & $\mathcal{M}_{h,\mu}^{(max)}$& $\alpha_{h,\mu}$ &
    $\mathcal{M}_{h,\mu}^{(max)}$& $\alpha_{h,\mu}$  & $\mathcal{M}_{h,\mu}^{(max)}$& $\alpha_{h,\mu}$ \\  \hline
    128 &  1.061E-03 & 1.167      &   6.221E-04& 1.334   & 3.327E-04 & 1.507\\
  256 &  4.727E-04& 1.167      &   2.473E-04& 1.334   & 1.196E-04 & 1.503\\
  512 &  2.106E-04 & 1.167 & 9.828E-05& 1.333 &4.296E-05& 1.500\\
  1024 & 9.381E-05& 1.167   & 3.905E-05& 1.333  &1.543E-05 &1.498\\
  2048 & 4.179E-05 &1.167  & 1.551E-05& 1.333  &5.534E-06& 1.496\\
  \hline
  $c_{\mu}$ & \multicolumn{2}{c|}{0.3049} & \multicolumn{2}{c|}{0.4024} & \multicolumn{2}{c|}{0.4975}\\
  \hline
\end{tabular}
  \centering
  \caption{Decay of the mixed power function for $\mu
   \in (2,3]$}\label{o3}
\end{table}

\begin{table}
\begin{tabular}{| c|c |c |c|c |}\hline
& \multicolumn{2}{c|}{$\mu = 10/3$} &  \multicolumn{2}{c|}{$\mu =
11/3$}\\ \hline
   $h^{-1}$   & $\mathcal{M}_{h,\mu}^{(max)}$& $\alpha_{h,\mu}$ &
    $\mathcal{M}_{h,\mu}^{(max)}$& $\alpha_{h,\mu}$   \\  \hline
    128 &  2.032E-04 & 1.491   &   1.661E-04 & 1.497\\
  256 &  6.995E-05& 1.497   &   5.808E-05 & 1.499\\
  512 &  2.419E-05 & 1.501 & 2.039E-05& 1.500\\
  1024 & 8.402E-06& 1.503 & 7.182E-06&1.501\\
  2048 & 2.919E-06 &1.505  & 2.523E-06& 1.502\\
  \hline
  $c_{\mu}$ & \multicolumn{2}{c|}{0.2814} & \multicolumn{2}{c|}{0.2368} \\
  \hline
\end{tabular}
  \centering
  \caption{Decay of the mixed power function for $\mu
   \in (3,4)$}\label{o4}
\end{table}

\theo{Conjecture}{conj}{ The mixed power function satisfies the
estimate {\rm{(\ref{ao})}} with the algebraic decay rate
\[
\alpha_{\mu}=\begin{cases}
\frac{\mu}{2}, \,\, & \,\, \textrm{for} \,\,\,\,  \mu\in (0,3)\setminus \{2\},\\
\frac{3}{2}, \,\, & \,\,  \textrm{for}\,\,\,\,\mu\in [3,4).
\end{cases}
\]
} 

\section{The mixed power function for $\mu = 3$}

Note that if Conjecture~\ref{conj} can be established for a
particular value $\mu\in [3,4)$, then the mixed power function bound
(\ref{mixed-bound}) implies a new and improved error estimate for
thin plate spline interpolation on the unit interval, namely, for
any $f \in C^{3}[0,1]$,
 \EQ{new-estimate}{ \max_{x\in\left[ 0,1\right]
}\left\vert f ( x ) -s_{h,2} ( x ) \right\vert =O\left(
h^{3/2}\right) ,\quad\text{as }h\rightarrow 0.}
This would provide a theoretical explanation of the numerical
results reported in \cite{B99,hm,Pow}.

In this section, we investigate Conjecture~\ref{conj} for the
special case $\mu =3$. By (\ref{mixed-pf}) and (\ref{mixedpfdir}),
the square of the mixed power function $\mathcal{M}_{h,3}$,
combining TPS Lagrange functions with the cubic basis function
$\phi_{3}$, is given by
\EQ{fourpow}{
\begin{aligned}
\mathcal{M}_{h,3}^{2}(x) & = \frac{A_{3}}{2\pi
}\int_{\R}\frac{\left\vert \Theta_{x,2}\left( t\right) \right\vert
^{2}}{\left\vert t\right\vert ^{4 }}dt
\\
& = \sum_{j=0}^{n}\sum_{k=0}^{n}\ell_{j,h}^{(2)}\left(
x\right)\ell_{k,h}^{(2)}\left(
x\right)|hj-hk|^{3}-2\sum_{j=0}^{n}\ell_{j,h}^{(2)}\left(
x\right)|x-hj|^{3}.
\end{aligned}
}

\begin{figure}[h]
\centering
\includegraphics[width=1.\textwidth]{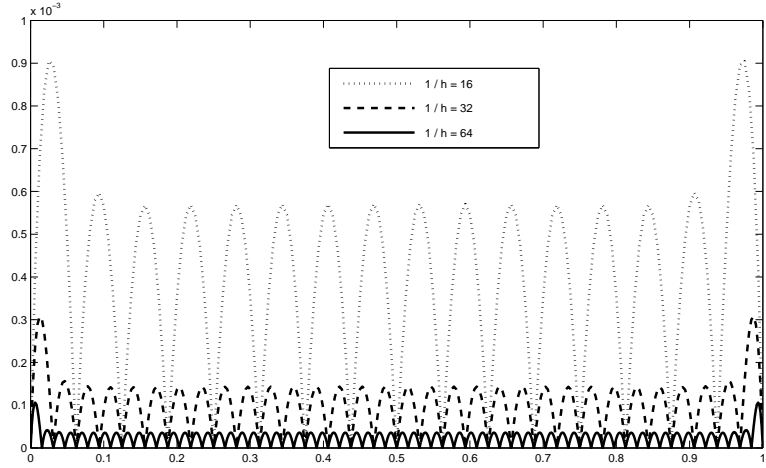}
\caption{Plot of $\mathcal{M}_{h,3}\left( x\right)$ for $h^{-1} =
16, 32$ and $64.$}\label{all}
\end{figure}

Figure \ref{all} illustrates the decay of $\mathcal{M}_{h,3}$ as
$h\rightarrow 0$. It can be confirmed numerically that the decay
rate $O(h^{3/2})$ of $\mathcal{M}_{h,3}$ suggested by Table \ref{o3}
applies uniformly on $[0,1]$, i.e.\ all peaks of the plot decay at
this rate.

We now relate the mixed power function $\mathcal{M}_{h,3}$ to a
classical error analysis method, namely the Peano kernel
representation. Let
\[
E_{h,x} (f) := f(x) - s_{h,2}(x) = f(x) - \sum\limits_{j=0}^{n}
f(hj) \ell_{j,h}^{(2)} (x).
\]
For each $x\in [0,1]$, $E_{h,x}$ is a continuous linear functional
on $C[0,1]$ with the usual max norm, and (\ref{replin}) implies that
the linear polynomials are in the null space of $E_{h,x}$. Then, for
any $f$ with an absolutely continuous first derivative on $[0,1]$,
Peano's theorem \cite[p.\ 271]{P81} implies
\EQ{peex}{f(x) - s_{h,2}(x) =
\int_{0}^{1}K_{h,x}(u)f^{\prime\prime}( u)\, du,\quad  \forall x \in
[0,1],}
where $K_{h,x}$ is the `Peano kernel' given by
\[
K_{h,x} (u) := (x-u)_{+} - \sum\limits_{j=0}^{n} \ell_{j,h}^{(2)}
(x) (hj-u)_{+}, \quad  u \in \R.
\]
\theo{Proposition}{peano}{
For each $ x \in [0,1]$, the mixed power function value
$\mathcal{M}_{h,3}\left( x \right)$ is a constant multiple of the
$L^{2} [0,1]$-norm of the Peano kernel $K_{h,x}$.
 }
\begin{proof}
The reproduction property (\ref{replin}) implies that $K_{h,x}$ is
compactly supported on $[0,1]$, and that
\[
K_{h,x} (u) = \frac{1}{2}\left(|x-u| - \sum\limits_{j=0}^{n}
\ell_{j,h}^{(2)} (x) |hj-u|\right), \quad  u \in \R.
\]
Then, by Lemma \ref{ftrbf}, the Fourier transform of the kernel
$K_{h,x}$ is the analytic and square integrable function
\[
\widehat{K_{h,x}}(t) = \frac{A_{1}}{2} \frac{\Theta_{x,2}\left( -t
\right)}{t^{2}}, \quad t\in\R,
\]
where $\Theta_{x,2}$ is defined by (\ref{theta}). Therefore, using
the first line of (\ref{fourpow}), the Parseval-Plancherel formula,
and the compact support of $K_{h,x}$, we deduce
\[
\begin{aligned}
\frac{2\pi}{A_{3}} \mathcal{M}_{h,3}^{2}\left(  x\right)
& = \int_{\R} \frac{\left\vert \Theta_{x,2}\left( -t \right)  \right\vert ^{2}}{t^{4}}\,dt \\
& = \frac{4}{A_1^2} \int_{\R}|\widehat{K_{h,x}}(t)|^{2}\,dt \\
& = \frac{8\pi}{A_1^2} \int_{0}^{1}|K_{h,x}(u)|^{2}\,du,
\end{aligned}
\]
which is the required conclusion.
\end{proof}

As a consequence, we obtain an alternative way of bounding the error
$f - s_{h,2}$ in terms of the mixed power function
$\mathcal{M}_{h,3}$, by using Cauchy-Schwarz directly in the
right-hand side of the Peano formula (\ref{peex}). The resulting
bound applies to any $f$ with an absolutely continuous first
derivative on $[0,1]$ and $f^{\prime\prime} \in L^2 [0,1]$. This
represents an improvement over
(\ref{mixed-bound})--(\ref{mixed-const}), which required $f\in C^3
[0,1]$ for $\mu=3$.

\begin{figure}[h]
\centering
\includegraphics[width=1.0\textwidth]{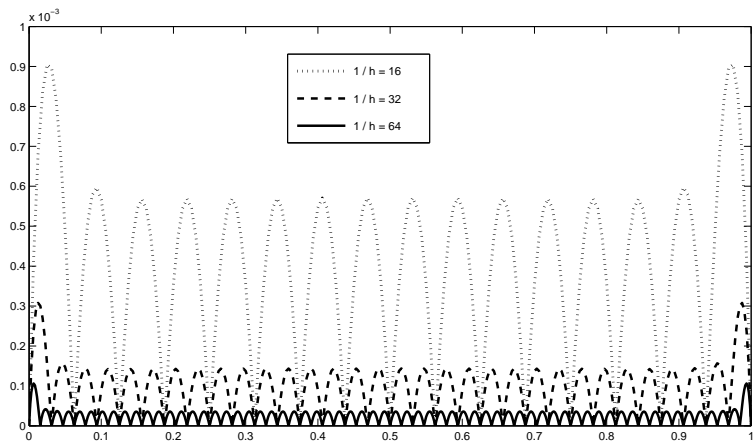}
\caption{Plot of $\mathcal{B}_{h}(x)$ for $h^{-1} = 16, 32$ and
$64.$}\label{allpea}
\end{figure}

Finally, a related question of interest is whether a sharper uniform
error bound can be obtained from (\ref{peex}) via H\"older's
inequality
\[
|f(x) - s_{h,2}(x)| \le \mathcal{B}_{h}(x) \left\Vert
f^{\prime\prime} \right\Vert_{ L^{\infty} [0,1] }, \quad x\in [0,1],
\]
where $f^{\prime\prime} \in L^{\infty} [0,1]$ and
\[
\mathcal{B}_{h}(x):=\int_{0}^{1}|K_{x,h}(u)|du.
\]
We note that this technique  was  used by Atkinson \cite{Atkins} in
the late 1960s
 to investigate the error behavior of
 natural cubic spline interpolant near the endpoints of the unit
 interval; see also Schaback \cite{Schabcub} for a treatment that is closer to
  our presentation. In the case of the TPS interpolant, a numerical
 answer to the question
is provided in Table  \ref{unibndB}, whose entries satisfy
\[
\mathcal{B}_{h} \left( \frac{h}{2} \right) =0.05059\, h^{\beta_{h}},\ \ \text{and }\ %
\mathcal{B}_{h} \left( \frac{1-h}{2} \right) =0.14955\,
h^{\sigma_{h}},
\]
i.e.,\ $\mathcal{B}_{h}$ decays approximately with the rate
$O(h^{3/2})$ near the endpoints of $[0,1]$ and this rate improves to
$O(h^2)$ near the midpoint. Also, Figure~\ref{allpea} shows that the
extreme peak value is well approximated by $\mathcal{B}_{h} \left(
\frac{h}{2} \right)$, while all of the lower peaks decay at the
faster rate. Therefore estimating the $L^1$-norm of the Peano kernel
leads to the same rate of decay $O(h^{3/2})$ of the uniform error
(\ref{uner})  as that predicted in (\ref{new-estimate}) by the mixed
power function $\mathcal{M}_{h,3}$.

We conclude the paper by remarking that any theoretical proof of the
uniform decay rate $O(h^{3/2})$ of $\mathcal{M}_{h,3} (x)$ or
$\mathcal{B}_{h} (x)$ for $x\in [0,1]$ will have to rely on specific
properties of the TPS Lagrange functions $\ell_{j,h}^{(2)}$,
$j\in\{0,1,\ldots,n\}$. A potentially useful such property is the
special case of \cite[Theorem~3.1]{ab01} stating that the
Lebesgue-type constant
\[
\max_{x\in [0,1]} \sum_{j=0}^{n} \left[ \ell_{j,h}^{(2)} (x)
\right]^{2}
\]
admits an upper bound independently of the mesh-size $h$. It remains
an open question whether this or other properties of the TPS Lagrange
functions can lead to further progress on the above conjectures.

\begin{table}
\begin{tabular}{| c|c |c |c |c|}
  \hline
   &  &  & &\\
   $h^{-1}$   & $\mathcal{B}_{h}\left(\frac{h}{2}\right)$ & $\beta_{h}$ &
  $\mathcal{B}_{h}\left(\frac{1-h}{2}\right)$ & $\sigma_{h}$ \\
  \hline
  64 &  1.024E-04 & 1.491 &3.633E-05& 2.001 \\
  128 &  3.533E-05 & 1.498 &9.098E-06& 2.001 \\
  256 &1.228E-05 &1.501  &2.293E-06& 1.999\\
  512 &  4.289E-06& 1.503 &5.694E-07& 2.000\\
  1024 & 1.502E-06 & 1.504 &1.434E-07& 1.999\\
  \hline
\end{tabular}
  \centering
  \caption{Decay of the $L^1$-norm of the Peano kernel.}\label{unibndB}
\end{table}

{\bf{Acknowledgements.}} A.B. is grateful to Michael Johnson (Kuwait
University) for a suggestion which improved the presentation of the
numerical results and for spotting a fallacious argument in a
tentative proof of Conjecture~\ref{conj} for $\mu =3$. The final
version has also benefited from the referees' suggestions.

                        \vskip1cm

                        %
%
                        %
%
%
%
%
%
%
%
%
                        %
%
                        %
%
%
                        %
%

                        {\large \bf Addresses:}\\[0.3cm]

            Aurelian Bejancu \\
                        Department of  Mathematics\\
                        Kuwait University \\
                        PO Box 5969\\
            Safat 13060 \\
                        Kuwait \\

                        Simon Hubbert (Corresponding Author) \\
                        School of Economics, Mathematics and Statistics\\
                        Birkbeck, University of London \\
                        Malet Street \\
                        London, WC1E 7HX\\
                        England \\

                        \end{document}